\documentclass[12pt]{article}
\usepackage{graphicx}
\usepackage{amsmath}
\usepackage{amsfonts}
\usepackage{amsthm}
\usepackage{esint}
\usepackage[T1]{fontenc}
\usepackage{url}
\usepackage{color}
\usepackage[margin=1in]{geometry}
\usepackage{amssymb,bm}
\usepackage{amsmath}
\usepackage{amsthm}
\usepackage{graphicx}
\usepackage[active]{srcltx} 
\usepackage{hyperref}
\hypersetup{pdfborder=0 0 0}

\newtheorem{theorem}{{\sc Theorem}}[section]

\newtheorem{remark}[theorem]{Remark}

\newtheorem{definition}[theorem]{Definition}

\def\XXint#1#2#3{{\setbox0=\hbox{$#1{#2#3}{\int}$ }
\vcenter{\hbox{$#2#3$ }}\kern-.6\wd0}}

\newcommand{\Ga}{\alpha}

\newcommand{\Ge}{\epsilon}

\newcommand{\Gg}{\gamma}

\newcommand{\Go}{\omega}

\newcommand{\GO}{\Omega}

\bmdefine\BGa{\alpha}
\bmdefine\BGb{\beta}
\bmdefine\BGd{\delta}
\bmdefine\BGe{\epsilon}
\bmdefine\BGve{\varepsilon}
\bmdefine\BGf{\phi}
\bmdefine\BGvf{\varphi}
\bmdefine\BGg{\gamma}
\bmdefine\BGc{\chi}
\bmdefine\BGi{\iota}
\bmdefine\BGk{\kappa}
\bmdefine\BGl{\lambda}
\bmdefine\BGn{\eta}
\bmdefine\BGm{\mu}
\bmdefine\BGv{\nu}
\bmdefine\BGp{\pi}
\bmdefine\BGth{\theta}
\bmdefine\BGvth{\vartheta}
\bmdefine\BGr{\rho}
\bmdefine\BGvr{\varrho}
\bmdefine\BGs{\sigma}
\bmdefine\BGvs{\varsigma}
\bmdefine\BGt{\tau}
\bmdefine\BGj{\tau}
\bmdefine\BGu{\upsilon}
\bmdefine\BGo{\omega}
\bmdefine\BGx{\xi}
\bmdefine\BGy{\psi}
\bmdefine\BGz{\zeta}
\bmdefine\BGD{\Delta}
\bmdefine\BGF{\Phi}
\bmdefine\BGG{\Gamma}
\bmdefine\BGL{\Lambda}
\bmdefine\BGP{\Pi}
\bmdefine\BGT{\Theta}
\bmdefine\BGS{\Sigma}
\bmdefine\BGU{\Upsilon}
\bmdefine\BGO{\Omega}
\bmdefine\BGX{\Xi}
\bmdefine\BGY{\Psi}

\bmdefine\BCA{{\mathcal A}}
\bmdefine\BCB{{\mathcal B}}
\bmdefine\BCC{{\mathcal C}}
\bmdefine\BCD{{\mathcal D}}
\bmdefine\BCE{{\mathcal E}}
\bmdefine\BCF{{\mathcal F}}
\bmdefine\BCG{{\mathcal G}}
\bmdefine\BCH{{\mathcal H}}
\bmdefine\BCI{{\mathcal I}}
\bmdefine\BCJ{{\mathcal J}}
\bmdefine\BCK{{\mathcal K}}
\bmdefine\BCL{{\mathcal L}}
\bmdefine\BCM{{\mathcal M}}
\bmdefine\BCN{{\mathcal N}}
\bmdefine\BCO{{\mathcal O}}
\bmdefine\BCP{{\mathcal P}}
\bmdefine\BCQ{{\mathcal Q}}
\bmdefine\BCR{{\mathcal R}}
\bmdefine\BCS{{\mathcal S}}
\bmdefine\BCT{{\mathcal T}}
\bmdefine\BCU{{\mathcal U}}
\bmdefine\BCV{{\mathcal V}}
\bmdefine\BCW{{\mathcal W}}
\bmdefine\BCX{{\mathcal X}}
\bmdefine\BCY{{\mathcal Y}}
\bmdefine\BCZ{{\mathcal Z}}

\bmdefine\Bzr{ 0}
\bmdefine\Ba{ a}
\bmdefine\Bb{ b}
\bmdefine\Bc{ c}
\bmdefine\Bd{ d}
\bmdefine\Be{ e}
\bmdefine\Bf{ f}
\bmdefine\Bg{ g}
\bmdefine\Bh{ h}
\bmdefine\Bi{ i}
\bmdefine\Bj{ j}
\bmdefine\Bk{ k}
\bmdefine\Bl{ l}
\bmdefine\Bm{ m}
\bmdefine\Bn{ n}
\bmdefine\Bo{ o}
\bmdefine\Bp{ p}
\bmdefine\Bq{ q}
\bmdefine\Br{ r}
\bmdefine\Bs{ s}
\bmdefine\Bt{ t}
\bmdefine\Bu{ u}
\bmdefine\Bv{ v}
\bmdefine\Bw{ w}
\bmdefine\Bx{ x}
\bmdefine\By{ y}
\bmdefine\Bz{ z}
\bmdefine\BA{ A}
\bmdefine\BB{ B}
\bmdefine\BC{ C}
\bmdefine\BD{ D}
\bmdefine\BE{ E}
\bmdefine\BF{ F}
\bmdefine\BG{ G}
\bmdefine\BH{ H}
\bmdefine\BI{ I}
\bmdefine\BJ{ J}
\bmdefine\BK{ K}
\bmdefine\BL{ L}
\bmdefine\BM{ M}
\bmdefine\BN{ N}
\bmdefine\BO{ O}
\bmdefine\BP{ P}
\bmdefine\BQ{ Q}
\bmdefine\BR{ R}
\bmdefine\BS{ S}
\bmdefine\BT{ T}
\bmdefine\BU{ U}
\bmdefine\BV{ V}
\bmdefine\BW{ W}
\bmdefine\BX{ X}
\bmdefine\BY{ Y}
\bmdefine\BZ{ Z}

\begin{document}
\title{On nonlocal geometric rigidity estimates in bounded uniform domains}
\author{Davit Harutyunyan\thanks{University of California Santa Barbara, harutyunyan@math.ucsb.edu} and 
Hayk Mikayelyan\footnote{University of Nottingham Ningbo China, Hayk.Mikayelyan@nottingham.edu.cn}}
\date{}

\maketitle

\begin{abstract}

We prove fractional geometric rigidity estimates and fractional Korn's first inequalities with general singular weights in bounded uniform domains in $\mathbb R^n.$ The classical choice of the energy well $SO(n)$ for the geometric rigidity estimate fails due to the lack of a notion of orientation for Lebesgue fields. The remedy is to choose $O(n)$ as the energy well. Due to the recently developed theory of local-to-nonlocal passage of rigidity estimates developed by the first author, both estimates reduce to the analogous local estimates on the cube, which we prove in several steps by compactness and linearization arguments. 
\end{abstract}

\section{Introduction}
\label{sec:1}

The Geometric rigidity expresses a fundamental principle of elasticity: a deformation whose nonlinear strain is small, must 
be close to a single rigid motion. The celebrated Geometric Rigidity Estimate (GRE) of Friesecke, James, and M\"uller [\ref{bib:Fri.Jam.Mue.}] quantifies this principle as follows.
  
\begin{theorem}[{[\ref{bib:Fri.Jam.Mue.}, Theorem~3.1]}]
Let $n\geq 2,$ let $p\in (1,\infty),$ and let $\Omega\subset\mathbb R^n$ be an open, bounded, connected, Lipschitz set. 
There exists a constant $C=C(n,p,\GO),$ such that for any vector field $ u\in W^{1,p}(\Omega,\mathbb R^n),$ there exists a 
proper rotation $\BR\in \mathrm{SO}(n)$ such that 
\begin{equation}
\label{1.1}
\|\nabla u-\BR \|_{L^p(\Omega)}^p \leq C\int_\GO \mathrm{dist}^p(\nabla u(x), \mathrm{SO}(n))dx.
\end{equation}
\end{theorem}

Since the work of Friesecke, James, and M\"uller, different variations, such as its extensions to mixed growth, incompatible fields, multiple wells, and variable domains illustrate the continuing importance of quantitative rigidity in variational problems; see 
[\ref{bib:Che.Con.},\ref{bib:DeL.Sze.},\ref{bib:Con.Dol.Mue.},\ref{bib:Mue.Sca.Zep.},\ref{bib:Fri.Kre.Zem.},\ref{bib:Har.Rod.}]. This estimate underlies the derivation of reduced elastic theories and provides a nonlinear counterpart of Korn's inequality. 
In these works, the corresponding local versions of the GRE have been proved for bounded connected Lipschitz domains. However, the original formulation (\ref{1.1}) of the GRE also holds in John domains, as the Whitney cube patching technique utilized in
[\ref{bib:Fri.Jam.Mue.},\ref{bib:Con.Zwi.}] goes through in such domains, due to the fact, that the weighted Poincar\'e inequality,
 with a weight that is the distance to the domain boundary raised to a power, holds in John domains, e.g., 
 [\ref{bib:Hurri},\ref{bib:Dre.Dur.}]. We emphasize though, that an explicit proof of this fact does not appear in the published literature, 
 yet the same proof as for instance in the paper [\ref{bib:Con.Zwi.}] for the so-called $(R,L)-$Lipschitz domains goes through for the GRE part.  At the same time, the validity, and even the formulation of the inequality in the nonlocal (fractional) setting is open in the contrary, while on the other hand the linear nonlocal (fractional) counterpart (Korn's inequality) has been settled in $C^1$ domains in [\ref{bib:Har.Men.Mik.Sco.}] and for bounded uniform domains in [\ref{bib:Aco.Dre.Dur.Lop.Oje.}], both for the standard kernel $\rho(x,y)=|x-y|^{-n-ps},$ where $s\in (0,1)$ and $p\in (1,\infty);$ see (\ref{1.4}) and also the papers  
[\ref{bib:Mengesha.1},\ref{bib:Men.Sco.1},\ref{bib:Rutkowski.1},\ref{bib:Har.Mik.}]. In this paper we formulate and prove the nonlocal GRE in bounded uniform domains with a general class of singular kernels, that also include the standard kernel
$\rho(x,y)=|x-y|^{-n-ps};$ the specifics are provided in Section~2 under main results. 

For nonlocal models, the deformation gradient is replaced by differences of the deformation at pairs of points. 
Recall, that for $s\in (0,1),$ and for a vector field $u\in L^p(\Omega,\mathbb R^n),$ the Gagliardo seminorm associated 
to the fractional $s-$derivative of the deformation $u$ is given by 
\begin{equation}
\label{1.2}
[u]_{W^{s,p}(\Omega)}^p=\int_{\Omega}\int_{\Omega}\frac{|u(x)-u(y)|^p}{|x-y|^{n+ps}}dx dy, 
\end{equation}
while the quantity 
$\| u\|_{W^{s,p}(\Omega)}=[ u]_{W^{s,p}(\Omega)}+\|u\|_{L^p(\Omega)}$ defines a norm 
[\ref{bib:Bou.Bre.Mir.1},\ref{bib:Maz.Sha.},\ref{bib:DiN.Pat.Val.}]. The fractional analogue of the linear elastic energy is given by (see [\ref{bib:Mengesha.1}]): 
\begin{equation}
\label{1.3}
[u]_{\mathcal{X}^{s,p}(\Omega)}^p=\int_{\Omega}\int_{\Omega}\frac{|(u(x)-u(y))\cdot(x-y)|^p}{|x-y|^{n+ps+p}}dxdy.
\end{equation}

The expression in (\ref{1.3}) has arisen in the theory of linear peridynamics [\ref{bib:Silling.1},\ref{bib:Silling.2},\ref{bib:Sil.Ept.Wec.Xu.Ask.},\ref{bib:Du.Zho.}] as the energy in the small strain regime in the analogy of linear elasticity, 
where the symmetric gradient $e(u)$ is the linear strain, and the integral $\int_\Omega |e(u)|^p dx$ is a replacement 
of to the linear elastic energy. It is not difficult to show that the condition $[u]_{\mathcal{X}^{s,p}(\Omega)}=0,$ which is 
equivalent to 
$$
(u(x)-u(y))\cdot (x-y)=0,\qquad\text{for a.e.}\qquad x,y\in\GO,
$$
implies even for $u\colon \GO\to\mathbb R^n$ fields with minimal regularity, like in the local case, that $u(x)=\BA x+b$ a.e. in $\GO$ for some skey-symmetric matrix $\BA\in\mathbb R^{n\times n}$ ($\BA+\BA^T=0$) and a vector $b\in\mathbb R^n; $ see [\ref{bib:Tem.Mir.}, Proposition~1.2] for an elementary proof. Fractional Korn's first inequality without boundary conditions, quantifies the above exact zero energy principle through the rigidity estimate
\begin{equation}
\label{1.4}
 \inf_{\Ba\in \mathfrak{se}(n)}[u-a]_{W^{s,p}(\Omega)}^p\leq C(\GO,p,n)[u]_{\mathcal{X}^{s,p}(\Omega)}^p
 \qquad\text{for all}\qquad u\in L^p(\Omega).
\end{equation}
Here, 
\begin{equation}
\label{1.5}
\mathfrak{se}(n)=\{a(x)=\BA x+b\colon\mathbb R^n\to\mathbb R^n\}
\quad\text{with}\quad \BA\in\mathbb R^{n\times n},\ \BA+\BA^T=0,\ b\in\mathbb R^n,
\end{equation}
is the space of infinitesimal rigid motions of $\mathbb R^n,$ see [\ref{bib:Har.Men.Mik.Sco.},\ref{bib:Aco.Dre.Dur.Lop.Oje.}]. 

It is known that measurable a.e. exact isometries from an open connected subset of $\mathbb R^n$ into $\mathbb R^n$ are 
of the form $r(x)=\BR x+b$, where $\BR\in \mathrm{O}(n),$ see [\ref{bib:Zal.Gat.}, Theorem~4.1]. This raises a corresponding 
rigidity question: does small distortion of pairwise distances force a deformation to be close to a single Euclidean isometry in a fractional Sobolev seminorm? The nonlinear strain must be chosen consistently with the interaction energy, and its zero set must 
be the space of exact isometries from $\GO$ into $\mathbb R^n.$ Consequently, one could naturally expect that an analogous to (\ref{1.4}) inequality should hold with the correct choice of the right-hand-side, that vanishes exactly on the space of exact 
isometries of $\mathbb R^n.$ The purpose of this paper is to establish a fractional geometric rigidity estimate based on 
bond-length distortion, including the endpoint $p=1.$ Recall the notation 
\begin{equation}
\label{1.6}
\mathfrak{so}(n)=\{r(x)=\BR x+b\colon\mathbb R^n\to\mathbb R^n\}
\qquad\text{with}\qquad \BR\in\mathrm{O}(n),\ b\in\mathbb R^n,
\end{equation} 
for the space of exact isometries of $\mathbb R^n.$ Hence, at first sight, a natural and naive choice for an appropriate fractional 
analogue of (\ref{1.1}) would be 
\begin{equation}
\label{1.7}
 \inf_{r \in \mathfrak{so}(n)}[u-r]_{W^{s,p}(\Omega)}^p \leq C(\GO,p,n)\int_{\Omega}\int_{\Omega}
 \frac{\bigl||u(x)-u(y)|-|x-y|\bigr|^p}{|x-y|^{n+ps}}dxdy, 
\end{equation}
for all $\Bu\in L^p(\Omega).$ The numerator of the nonlinear energy measures the change in the length of the bond joining $x$ and $y$. Unlike a gradient-based strain, this expression is meaningful for measurable deformations and does not require a classical derivative. The use of bond deformation as a basic kinematic quantity connects (\ref{1.7}) with peridynamics. 
In the present manuscript, we will prove, that indeed (\ref{1.7}) holds in more general form and terms imposed on the kernel and the domain. Thus the estimate (\ref{1.7}) supplies a nonlinear coercivity statement for a particular singular bond kernel.  
The relationship with the fractional Korn inequality is revealed by linearization like in the local case by taking for a deformation 
the well-known field $u_\varepsilon(x)=x+\varepsilon v(x)$ in (\ref{1.7}) and sending $\varepsilon$ to zero. This will yield 
 (\ref{1.4}) from (\ref{1.7}) through a routine argument as in the classical case. However, we will derive (\ref{1.4}) from the local
 counterpart via the general theory in [\ref{bib:Harutyunyan}] described below, as that local Korn's inequality will be an important 
 component of the proof of the corresponding local geometric rigidity estimate via linearization. 
 
We will prove (\ref{1.7}) and  (\ref{1.4}) utilizing the theory of reduction of nonlocal general operator inequalities with singular kernels to the analogous local estimates on the cube. This framework was recently developed by the first author in [\ref{bib:Harutyunyan}].  

While we were at final stages of typing and polishing the present manuscript, a formulation and a proof of the fractional geometric rigidity estimate was open in the literature, and Tadele Mengesha and James Scott independently (from us) proved versions of the fractional Korn's inequality in [\ref{bib:Men.Sco.2}] and Scott formulated and proved the nonlocal geometric rigidity estimate in [\ref{bib:Scott}], in both papers in bounded Lipschitz domains, and kindly shared a pdf copy of their papers with us (and posted preprints on arxiv.org). Mengesha and Scott prove (\ref{1.4}) in Lipschitz domains with the peridynamics kernel with optimal constants scaling in terms of $1-s$ as $s\to 1-$ also in the endpoint cases $p=1,\infty;$ see [\ref{bib:Men.Sco.2}, Theorems~1.1-1.3]. Scott proves (\ref{1.7}) in bounded Lipschitz domains with more general radial kernels and for $1\leq p<\infty;$ see [\ref{bib:Scott}, Theorem~2.1]. We will prove both (\ref{1.7}) and (\ref{1.4}) with general radial kernels, for $1\leq p<\infty,$ and in the class of 
bounded uniform domains, which contains the class of bounded Lipschitz domains, see Section~2 for details. It is worth 
mentioning, that our approaches and proof schemes significantly (almost completely) differ from each other. 

In conclusion, note, that an essential distinction from (\ref{1.1}) concerns orientation. That is, if one takes $\mathrm{SO}(n)$ 
instead of $\mathrm{O}(n)$ in (\ref{1.7}), then it would fail, due to the possibility of choosing an isometry of the form 
 $r(x)=(\BI-2e_1\otimes e_1)\BR x,$ with $\BR\in \mathrm{SO}(n),$ which is an orientation altering orthogonal 
 transformation, that is en exact isometry of $\mathbb R^n.$ We discuss this in Remark \ref{Rem:3.2}.


\section{Problem setting and the main results}
\setcounter{equation}{0}
\label{sec:2}

In what follows $\GO\subset\mathbb R^n$ will be a bounded uniform domain. The class of uniform domains is a subclass of John 
domains, the definition of which is given below. 


\begin{definition}[Uniform domain] A domain $\GO\subset\mathbb R^n$ is called a uniform domain, if there exists a constant 
$C\geq 1,$ such that for every two points $x,y\in\GO,$ there exists a rectifiable curve $\Ga\colon[0,s]\to \GO,$ 
that starts at $x,$ ends at $y,$ and satisfies $\ell(\Gg)\leq C|x-y|,$ 
$\min(\ell(x,\Gg(t)),\ell(y,\Gg(t)))\leq C\mathrm{dist}(\Gg(t),\partial\GO)$ for all $t\in [0,s],$ where $\ell(z,\Gg(t))$ is 
the length of the curve between the points $z$ and $\Gg(t)$ for $z=x,y.$
\end{definition}

We refer the reader to the original paper of John [\ref{bib:John}] for details on John domains and to the papers 
[\ref{bib:Mar.Sar.},\ref{bib:Geh.Osg.},\ref{bib:Jones},\ref{bib:Seeger}] for more details on uniform domains. The class of Lipschitz 
domains is a proper subset of the class of uniform domains, while the later is a proper subset of the class of John domains. 
Also, uniform domains are open, bounded and simply connected by definition. 

As already mentioned, we will prove (\ref{1.7}) and (\ref{1.4}) with a class of general kernels $\rho(x,y),$ see Theorems~2.1 
and 2.2. For that purpose, we recall below the conditions imposed on the kernel in [\ref{bib:Harutyunyan}] under which the 
general localization theory is valid. We assume the kernel 
\begin{equation}
\label{2.1}
\rho(x,y)\colon\GO\times \GO\to (0,\infty)
\end{equation}
satisfies the hypotheses below. \\


\noindent\textbf{The hypotheses on the kernel $\rho$.} 
There exist constants $c_1>0$ and $0<c_2 \leq c_2'<1$ such that for all $x,y \in \mathbb R^n$ one has
\begin{itemize}
\item[(H1)] 
\begin{equation}
\label{2.2}
\rho(y) \leq c_1\rho(x), \qquad \text{if} \qquad |x|\leq |y|.
\end{equation}
\item[(H2)] 
\begin{equation}
\label{2.3}
c_2 \rho(x) \leq \rho(2y) \leq c_2' \rho(x),\qquad \text{if} \qquad |x|=|y|\neq 0.
\end{equation}
\end{itemize}
Now, by notational abuse, if we define for $t>0,$
\begin{equation}
\label{2.4}
\rho(t)=\rho(te_1),
\end{equation}
then by (H1)-(H2) we get an equivalent radial kernel that satisfies $\rho(t)\colon (0,\infty)\to (0,\infty)$ and we get the analogous one dimensional conditions:
\begin{itemize}
\item[(H1')]
\begin{equation}
\label{2.5}
\rho(t_2) \leq c_1\rho(t_1)\qquad \text{for all}\qquad 0<t_1\leq t_2<\infty. 
 \end{equation}
 \item[(H2')]
\begin{equation}
\label{2.6}
c_2\rho(t) \leq \rho(2t)\leq c_2' \rho(t)\qquad \text{for all}\qquad t>0.
 \end{equation}
 The moment conditions below are also assumed.
\item[(H3)] There exists a constant $c_3>0$ such that
\begin{equation}
\label{2.7}
\frac{1}{a^{n+p}\rho(a)}\int_0^a t^{n+p-1}\rho(t)dt\leq c_3 \qquad\text{for all}\qquad a>0.
\end{equation}
\item[(H$\rho$4)] There exists a constant $c_4>0$ such that 
\begin{equation}
\label{2.8}
\int_a^1 t^{n-1}\rho(t)dt\leq c_4a^n\rho(a)\qquad\text{for all}\qquad 0<a<1.
\end{equation}
\end{itemize}
The standard kernel satisfies the above hypotheses with the constants 
\begin{equation}
\label{2.9}
c_1=1,\quad c_2=c_2'=2^{-n-ps},\quad  c_3=\frac{1}{p(1-s)},\quad c_4=\frac{1}{ps}.
\end{equation}

For a field $u\in L^1_{loc}(\GO,\mathbb R^n),$ denote for the sake of brevity 
\begin{equation}
\label{2.10}
\Go(u)=\Go(u(x,y)):=|u(x)-u(y)|-|x-y|,\qquad x,y\in\GO.
\end{equation}

The first theorem is the nonlocal geometric rigidity estimate (GRE) on uniform domains. 

\vspace{0.2cm}

\begin{theorem}[Nonlocal GRE in uniform domains]
\label{Thm:2.2}
Let $n\geq 2$ and $p\in [1,\infty)$ and let the kernel $\rho$ satisfy the hypotheses (H1')-(H4). 
Then there exists a constant $C$ that depends on $n,p,c_1,c_2,c_2',c_3,c_4$ and $\GO,$ such that 
\begin{equation}
\label{2.11}
 \inf_{r \in \mathfrak{so}(n)}\int_\GO\int_\GO\rho(x,y) \bigl|u(x)-u(y)-[r(x)-r(y)]\bigr|^pdxdy
  \leq C \int_{\GO}\int_{\GO} \rho(x,y) |\Go(u(x,y))|^pdxdy, 
\end{equation}
for all vector fields $u\in L^1_{loc}(\Omega,\mathbb R^n).$ 
\end{theorem}

Next is the Korn's first inequality (KFI) without boundary conditions on uniform domains. 

\vspace{0.2cm}

\begin{theorem}[Nonlocal KFI in uniform domains]
\label{Thm:2.3}
Let $n\geq 2$ and $p\in [1,\infty)$ and let the kernel $\rho$ satisfy the hypotheses (H1')-(H4). 
Then there exists a constant $C$ that depends on $n,p,c_1,c_2,c_2',c_3,c_4$ and $\GO,$ such that 
\begin{align}
\label{2.12}
 \inf_{a \in \mathfrak{se}(n)}\int_\GO\int_\GO\rho(x,y)& \bigl|u(x)-u(y)-[a(x)-a(y)]\bigr|^pdxdy \\ \nonumber
 & \leq C \int_{\GO}\int_{\GO} \rho(x,y) \frac{\bigl|(u(x)-u(y))\cdot(x-y)\bigr|^p}{|x-y|^p}dxdy, 
\end{align}
for all vector fields $u\in L^1_{loc}(\Omega,\mathbb R^n).$ 
\end{theorem}

Next we recall the general singular kernel removal theory in [\ref{bib:Harutyunyan}, Theorems~2.1 and 2.3]. The following is a combination of  Theorems~2.1 and 2.3 in [\ref{bib:Harutyunyan}], see also the hypotheses (HA1)-(HA2) on the operator $\BA$ and discussion on the operators $\BA$ and $\BB$ and the fractional Korn and geometric rigidity inequalities in Section~4 of [\ref{bib:Harutyunyan}].

\vspace{0.2cm}

\begin{theorem}[{[\ref{bib:Harutyunyan}, Theorems~2.1 and 2.3]}]
\label{Thm:2.4}
Let $n\geq 2$ and $1\leq p<\infty,$ let $\GO\subset\mathbb R^n$ be a bounded uniform domain. Assume the kernel 
$\rho$ satisfies hypotheses (H1')-(H4). Let $\mathcal{F}\subset\mathcal{P}(\GO)$ be the set of all open and connected subsets of 
$\GO.$ Assume the pair of operators $\BA\colon\mathcal{F}\times L^1_{loc}(\GO,\mathbb R^n)\to L^1_{loc}(\GO,\mathbb R^{n})$ and $\BB \colon L^1_{loc}(\GO,\mathbb R^n)\to L^1_{loc}(\GO \times \GO,\mathbb R^{n})$ is 
defined below as in (i), for the fractional Korn inequality, or in (ii), for the fractional geometric rigidity estimate 
\begin{itemize} 
\item[(i)] Define 
\begin{equation}
\label{2.13}
\BA(E,u)(x)=u(x)-a_E(x),\quad\text{and}\quad (\BB u)(x,y)=\frac{(u(x)-u(y))\cdot(x-y)}{|x-y|} \quad x,y\in\GO,
\end{equation}
where $a_E(x)$ is the minimizer infinitesimal rigid motion in (\ref{2.18}) within the set $E,$ that exists as $\mathfrak{se}(n)$ 
is finite dimensional. 
\item[(ii)]  Define 
\begin{equation}
\label{2.14}
\BA(E, u)(x)=u(x)-r_E(x),\quad\text{and}\quad (\BB u)(x,y)=\Go(u(x,y)), \quad x,y\in\GO,
\end{equation}
where $r_E(x)$ is the minimizer infinitesimal rigid motion in (\ref{2.17}) within the set $E,$ that exists as $\mathfrak{so}(n)$
 is finite dimensional. 
\end{itemize}
Then statement (iii) implies statement (iv) below. 
\begin{itemize}
\item[(iii)] There exists a constant $C=C(n,p)$ such that for every cube $Q\subseteq \GO$ with $4Q\subseteq \GO,$
the unweighted inequality 
\begin{equation}
\label{2.15}
\int_Q |A(Q,u)(x)|^pdx \leq \frac{C(n,p)}{|Q|} \int_{4Q}\int_{4Q} |(\BB u)(x,y)|^pdxdy, 
\end{equation}
for all vector fields $u\in L^1_{loc}(4Q,\mathbb R^n).$
\item[(iv)] There exists a constant $\widetilde C=\widetilde  C(n,p,\rho,\GO)$ and a non-degenerate cube $Q\subset \GO,$ such that 
the inequality
\begin{align}
\label{2.16}
\int_{\Omega}\int_{\Omega}\rho(x,y)\bigl|A(Q,u)(x)-A(Q,u)(y) \bigr|^pdxdy
\leq \widetilde  C\int_{\Omega}\int_{\Omega} \rho(x,y)|(\BB u)(x,y)|^pdxdy,
\end{align}
holds for all $u\in L^1_{loc}(\GO,\mathbb R^N).$ 
\end{itemize}

\end{theorem}

Consequently it is sufficient to prove the local unweighted estimates on nested concentric cubes (\ref{2.15}) for both definitions (i) and (ii) of the operators $\BA$ and $\BB.$ The theorems are formulated below.

\vspace{0.2cm}

\begin{theorem}[Local GRE on concentric cubes]
\label{Thm:2.5}
Let $n\geq 2$ and $p\in [1,\infty)$ and let $Q\subset\mathbb R^n$ be a cube. Then there exists a constant $C$ that depends 
on $n$ and $p$ such that 
\begin{equation}
\label{2.17}
 \inf_{r \in \mathfrak{so}(n)} \int_Q |u(x)-r(x)|^pdx \leq \frac{C(n,p)}{|Q|} \int_{4Q}\int_{4Q} |\Go(u(x,y))|^pdxdy, 
\end{equation}
for all vector fields $u\in L^1_{loc}(4Q,\mathbb R^n).$ 
\end{theorem}

\vspace{0.2cm}

\begin{theorem}[Local KFI on concentric cubes]
\label{Thm:2.6}
Let $n\geq 2$ and $p\in [1,\infty)$ and let $Q\subset\mathbb R^n$ be a cube. Then there exists a constant $C$ that depends 
on $n$ and $p$ such that
\begin{align}
\label{2.18}
 \inf_{a \in \mathfrak{se}(n)}\int_Q|u(x)-a(x)|^pdx \leq  \frac{C(n,p)}{|Q|} \int_{4Q}\int_{4Q} \frac{\bigl|(u(x)-u(y))\cdot(x-y)\bigr|^p}{|x-y|^p}dxdy, 
\end{align}
for all vector fields $u\in L^1_{loc}(4Q,\mathbb R^n).$ 
\end{theorem}


 \section{Proof of the local inequalities on nested cubes}
\setcounter{equation}{0}
\label{Sec:3}
 
 We first prove Theorem~2.6, as in the compactness argument for (\ref{2.17}) we will need (\ref{2.18}) via linearization. 
 
\subsection{Proof of Theorem~2.6}
\label{Sub:3.1}

\begin{proof}[Proof of Theorem~2.6] By translation and rescaling, we can assume without loss of generality that $Q$ is the centered unit cube $Q=(-\tfrac12,\tfrac12)^n$. Define for measurable fields $u,$ for the sake of brevity, the directional difference 
\begin{equation}
\label{3.1} 
E u(x,y)=(u(x)-u(y))\cdot\frac{x-y}{|x-y|},
\end{equation}
that is similar to the classical symmetrized gradient, i.e., we can call is a zeroth order strain. We will argue by contradiction and
assume that (\ref{2.18}) fails. Then, by homogeneity, and the fact that the minimum on the left of (\ref{2.18}) is attained, we can assume that for some $\{u_j\}\subset L^1_{loc}(4Q,\mathbb R^n)$ and $a_j\in\mathfrak{se}(n)$ the conditions hold:

\begin{equation}
\label{3.2}   
\lim_j\|Eu_j\|_{L^p(4Q\times 4Q)}=0,\qquad\text{and}
\qquad  \inf_{a\in\mathfrak{se}(n)} \|u_j-a\|_{L^p(Q)}=\|u_j-a_j\|_{L^p(Q)}=1. 
\end{equation}
Substituting $v_j=u_j-a_j$ we have the normalized conditions
\begin{equation}
\label{3.3}   
\lim_j\|Ev_j\|_{L^p(4Q\times 4Q)}=0,\qquad\text{and}
\qquad  \inf_{\Ba\in\mathfrak{se}(n)} \|v_j-a\|_{L^p(Q)}=\|v_j\|_{L^p(Q)}=1, 
\end{equation}
which imply that $v_j\in L^p(Q).$ Through an application of the generalized Riesz-Fr\'echet-Kolmogorov criterion, we will prove
that in fact (\ref{3.3}) implies that $\{v_j\}$ is strongly relatively compact in $L^p(Q).$ The boundedness of $\{v_j\}$ in $L^p(Q)$ 
 is in (\ref{3.3}), hence we only need to verify the so-called joint translational equi-continuity of $v_j,$ see [\ref{bib:Han.Hol.}, Remark~6] and (\ref{3.10}) below. Following the calculations in [\ref{bib:Tem.Mir.}, Proposition~1.2], we can write 
$$v(x)x-v(x)y-v(y)x+v(y)y=|x-y|Ev(x,y),$$
thus writing for the pair $(x,z)\in Q\times Q$ and subtracting we obtain 
\begin{equation}
\label{3.4}   
v(x)(z-y)=x(v(y)-v(z))+zv(z)-yv(y)+|x-y|Ev(x,y)-|x-z|Ev(x,z).
\end{equation}
Writing again (\ref{3.4}) for the triple $(x+h,y,z)$ and subtracting we get 
\begin{align}
\label{3.5}   
&(v(x+h)-v(x))\cdot (z-y)=h\cdot(v(y)-v(z))\\ \nonumber
&+|x+h-y|Ev(x+h,y)-|x+h-z|Ev(x+h,z)-|x-y|Ev(x,y)+|x-z|Ev(x,z).
\end{align}
Now for every triple $(x,y,h)$ with $x,y,x+h\in Q,$ denote
\begin{equation}
\label{3.6}   
U(x,y,h)=\left\{z\in Q \ : \ \left|(v(x+h)-v(x))\cdot (z-y)\right|\geq \frac{1}{10}|v(x+h)-v(x)|\right\}.
\end{equation}
Clearly $U(x,y,h)$ is measurable, and that there exists a dimensional constant $c_0(n)$>0, such that for every fixed triple $(x,y,h)$ with $x,x+h\in Q$ and $y\in\frac{1}{2}Q$ one has $|U(x,y,h)|\geq c_0(n).$ Now, fix such a triple $(x,y,h)$ and integrate (\ref{3.5}) 
in $z$ over $U(x,y,h)$ to get with an application on Minkowski's inequality on the right-hand side:
\begin{align}
\label{3.7}   
\frac{c_0(n)}{10^p} &|v(x+h)-v(x)|^p\leq C(p)|h|^p\int_Q |v(y)-v(z)|^pdz+C(n,p)|Ev(x+h,y)|^p\\ \nonumber
&+C(n,p)\int_Q |Ev(x+h,z)|^pdz+C(n,p)|Ev(x,y)|^p+C(n,p)\int_Q |Ev(x,z)|^pdz.
\end{align} 
Another integration in $(x,y)\in Q\times \frac{1}{2}Q$ such that $x+h\in Q,$ yields 
\begin{align}
\label{3.8}   
\frac{c_0(n)}{10^p\cdot 2^n}\int_Q&|v(x+h)-v(x)|^pdx \leq C(n,p)|h|^p\int_Q \int_Q |v(y)-v(z)|^pdzdy\\ \nonumber
&+C(n,p)\int_{Q} \int_{Q}|Ev(x+h,y)|^pdxdy +C(n,p)\int_{Q} \int_{Q}|Ev(x+h,z)|^pdzdx \\ \nonumber
&+C(n,p)\int_{Q} \int_{Q}|Ev(x,y)|^pdxdy+C(n,p) \int_{Q} \int_{Q}|Ev(x,z)|^pdzdx  \\ \nonumber
&\leq C(n,p)|h|^p\int_Q |v(y)|^pdy+C(n,p)\int_{Q} \int_{Q}|Ev(x,y)|^pdxdy.
\end{align}  
Applying (\ref{3.8}) to the field $v_j$ we get thanks to (\ref{3.3}):
\begin{equation}
\label{3.9}   
\int_Q |v_j(x+h)-v_j(x)|^pdx\leq C(n,p)|h|^p+C(n,p)\int_{Q} \int_{Q}|Ev_j(x,y)|^pdxdy,
\end{equation} 
which yields the desired joint translational equi-continuity 
\begin{equation}
\label{3.10}   
\int_Q|v_j(x+h)-v_j(x)|^pdx \longrightarrow 0\qquad\text{as} \quad h\to 0, \ j\to\infty.
\end{equation} 
 This is not yet sufficient for a direct application of the Fr\'echet-Kolmogorov criterion, but the proof of the later can be 
 adjusted to treat the joint limit situation (\ref{3.10}). Consequantly, by  [\ref{bib:Han.Hol.}, Remark~6], we infer that for a 
 (non-relabeled) subsequence $v_j\to v$ in $L^p(Q)$ for some $v\in L^p(Q),$ which yields by the triangle inequality that 
 we have for $v$ the zero energy condition:
 \begin{equation}
\label{3.11}   
\|Ev\|_{L^p(Q\times Q)}=0.
\end{equation}
As already mentioned before, the kernel of $E$ is the space of infinitesimal rigid motions; see [\ref{bib:Tem.Mir.}, Proposition~1.2] 
for details. Thus we conclude that $v\in \mathfrak{se}(n)$ on $Q,$ but in that case we would have due to (\ref{3.3}), that 

$$
1=\|v_j\|_{L^p(Q)}=\inf_{\Ba\in\mathfrak{se}(n)} \|v_j-a\|_{L^p(Q)}\leq \|v_j-v\|_{L^p(Q)}\to 0,
$$
which is a contradiction. This completes the proof of Theorem~2.6.

\end{proof}

\begin{remark}
\label{Rem:3.1}
In fact, we have proven the inequality (\ref{2.18}) without quadrupling the cubes, but rather over the same cubes $Q$ on the right, i.e., we proved the estimate  
\begin{align}
\label{3.12}
 \inf_{a \in \mathfrak{se}(n)}\int_Q|u(x)-a(x)|^pdx \leq  \frac{C}{|Q|} \int_{Q}\int_{Q} \frac{\bigl|(u(x)-u(y))\cdot(x-y)\bigr|^p}{|x-y|^p}dxdy, 
\end{align}
for all vector fields $u\in L^1_{loc}(Q,\mathbb R^n).$ 
\end{remark}


\subsection{Proof of Theorem~2.5}
\label{Sub:3.2}

\begin{proof}[Proof of Theorem~2.5]
As for the case of (\ref{2.18}), we will prove the stronger estimate without enlarging the cubes, i.e.,
 \begin{equation}
\label{3.13}
\inf_{r\in\mathfrak{so}(n)} \fint_Q |u(x)- r(x)|^p dx \le C(n,p) \fint_Q \fint_Q |\Go(u(x,y))|^p dx dy,
\end{equation}
for all vector fields $u\in L^1_{loc}(Q,\mathbb R^n),$ where again by a variable translation and rescaling we can assume that $Q$
is the centered unit cube. Assume again in contradiction that (\ref{3.13}) fails. 
This yields a sequence $\{u_j\}\subset L^p(Q,\mathbb R^n)$, which we translate to have zero means on $Q$ without affecting (\ref{3.13}) and with
\begin{equation}
\label{3.14}  
a_j=\min_{r\in\mathfrak{so}(n)}\|u_j-r(x)\|_{L^p(Q)}> j d_j,
 \qquad \text{where}\qquad d_j=\|\Go(u_j)\|_{L^p(Q\times Q)}.
\end{equation}
Note, that unlike (\ref{3.12}), inequality (\ref{3.13}) is not homogeneous in $u$ and we can't take apriorily $a_j=1$ and $d_j\to 0.$
We will show first that in fact (\ref{3.14}) has enough structure to force $d_j,a_j\to 0$ as $j\to\infty.$ The minimum in (\ref{3.14}) 
exists by compactness of $\mathrm{O}(n)$, $u_j\in L^p(Q)$ and coercivity in $b$. We can write for $(x,y)\in Q^2$ that 
$$
u_j(x)=\int_Q(u_j(x)-u_j(y))dy,
$$ 
thus we have by Jensen 
$$
|u_j(x)|^p\leq \int_Q|u_j(x)-u_j(y)|^pdy,
$$
which integrating in $x\in Q$ and applying the triangle inequality and Minkowski as
$$
|u_j(x)-u_j(y)|^p\leq 2^{p-1}|\Go(u_j)|^p+2^{p-1}|x-y|^p,
$$
which we can integrate in $x$ over a measurable set $A\subset Q$ to discover
\begin{align}
\label{3.15}
\int_A|u_j(x)|^pdx & \leq \int_A \int_Q|u_j(x)-u_j(y)|^pdydx\\ \nonumber
&\leq 2^{p-1}d_j+C(n,p)|A|.
\end{align}
We get from (\ref{3.15}) by the triangle inequality and taking $A=Q$ that
$$
 a_j\le\|u_j-x\|_{L^p(Q)}\le C(n,p)(1+d_j),
$$
thus for large enough indexes $j$, we get from $a_j> jd_j,$ the upper bound 
$$
d_j\le \frac{C(n,p)}{j-C(n,p)}\longrightarrow 0.
$$ 
Now we aim to prove again that $\{u_j\}$ is strongly pre-compact in $L^p(Q).$ The bound (\ref{3.15}) is not yet sufficient for 
strong compactness, thus we pass to the mollifications to gain access to embedding results. Let $\Ge>0$ be small and let 
$\eta_\Ge$ be the standard smooth radial mollifier supported on the ball $B_\Ge(0)$ and let 
$u_j^\Ge=u_j\ast \eta_\Ge\in C^\infty(Q_\Ge,\mathbb R^n)$ be the mollification of $u_j$ defined on the cube $Q$ minus
the $\Ge-$collar $Q_\Ge=\{x\in Q \  :  \  \mathrm{dist}(x,\partial Q)>\Ge\},$ see [\ref{bib:Eva.Gar.}, Section 4.2] for details.
We will first show that 
\begin{equation}
\label{3.16}
 \|u_j-u_j^\Ge\|_{L^p(Q_\Ge)}^p \leq C(n,p)(\Ge^p+\Ge^{-n}d_j^p).
\end{equation} 
 We have for every $x\in Q_\Ge$ and $z\in B_\Ge(0)$ that $x-z\in Q,$ thus we can write 
 $$
 u_j(x)-u_j^\Ge(x)=\int_{B_\Ge(0)}\eta_\Ge(z)(u_j(x)-u_j(x-z))dz.
 $$
 Since $\eta_\Ge(z)dz$ is a probability measure, we have by Jensen and Minkowski, that 
 \begin{align*}
  |u_j(x)-u_j^\Ge(x)|^p & \leq \int_{B_\Ge(0)}\eta_\Ge(z)|u_j(x)-u_j(x-z)|^pdz\\ 
 & \leq 2^{p-1}\int_{B_\Ge(0)}\eta_\Ge(z)\left(|z|^p+|\Go(u_j(x,x-z))|^p\right)dz,  
 \end{align*}
 which integrating over $Q_\Ge$ we arrive at the desired estimate:
  \begin{align*}
   \|u_j-u_j^\Ge\|_{L^p(Q_\Ge)}^p  & \leq 2^{p-1}  |Q_\Ge|\int_{B_\Ge(0)}\eta_\Ge(z)|z|^pdz+
  2^{p-1}\int_{Q_\Ge}\int_{B_\Ge(0)}\eta_\Ge(z)\Go(v_j(x,x-z))^pdzdx \\
  &\leq 2^{p-1}\Ge^p+2^{p-1}\|\eta\|_{\infty}\Ge^{-n}d_j^p.
  \end{align*}
  Now, for any fixed $\Ge>0,$ the sequence $u_j^\Ge$ is bounded in $W^{1,p}(Q_\Ge)$ (see [\ref{bib:Eva.Gar.}, Theorem~4.1]), namely one has
  $$
  \|u_j^\Ge\|_{W^{1,p}(Q_\Ge)}\leq C (n,p,\Ge)\|u_j\|_{L^p(Q)}. 
  $$
 Thus by the Rellich-Kondrachov theorem it is strongly precompact in $L^p(Q_\Ge).$ We now extend each $u_j^\Ge$ by zero 
 in the collar $Q\setminus Q_\Ge$ and keep the same notation for the extension. Then from (\ref{3.16}) and (\ref{3.15}) 
 we discover 
 \begin{equation}
\label{3.17}
 \limsup_{j\to\infty}\|u_j-u_j^\Ge\|_{L^p(Q)}^p \leq C(n,p)\Ge.
\end{equation}  
 Hence, we get the strong precompactness of $\{u_j\}$ in $L^p(Q)$ by choosing now $\Ge_k\to 0$ and applying a diagonal argument together with the joint equi-integrability condition in (\ref{3.15}). Consequently, for some $u\in L^p(Q)$ we have
\begin{equation}
\label{3.18}
 u_j(x)\longrightarrow u \qquad\text{strongly in}\qquad  L^p(Q). 
\end{equation}
Next we prove that $u=\BR x$ for some $\BR \in \mathrm{O}(n).$ The triangle inequality gives
$$
 \|\Go(u)-\Go(u_j)\|_{L^p(Q\times Q)} \le 2\|u-u_j\|_{L^p(Q)},
$$
thus we get from (\ref{3.14}) and (\ref{3.18}) that $\|\Go(u)\|_{L^p(Q\times Q)}=0,$ and hence 
\begin{equation}
\label{3.19}
|u(x)-u(y)|=|x-y|\qquad \text{ for almost every pair} \qquad x,y\in Q.
\end{equation}
Theorem~4.1 of [\ref{bib:Zal.Gat.}] with the choice of the probability measure $\mu=\chi_Qdx$ implies that $u(x)=\BR x+b$ for 
some $\BR\in\mathrm{O}(n)$ and $b\in \mathbb R^n.$ The fact that $Q$ is centered at the origin and the zero mean 
condition on $u$ give $b=0.$ This yields $u(x)=\BR x$ a.e. in $Q.$ Consequently, we also get in (\ref{3.14}) that 
$$a_j=\min_{r\in\mathfrak{so}(n)}\|u_j-r(x)\|_{L^p(Q)}\leq \|u_j-u(x)\|_{L^p(Q)}\to 0.$$
We will now finish the proof by linearization. If in (\ref{3.14}) we choose minimizing rotations $\BR_j$ and replace $u_j$ by $\BR_j^T(u_j-b_j),$ 
the values of $a_j$ and $d_j$ will not change. For this new sequence, the minimizing rotations reduce to the identity matrix, 
thus can normalize 
\begin{equation}
\label{3.20}
 u_j(x)=x+a_jv_j(x) \qquad\text{with}\qquad \|v_j\|_{L^p(Q)}=1.
\end{equation}
Due to the condition $a_j\to 0,$ we aim to linearize $\Go(u)=|x-y+a_j(v_j(x)-v_j(y)|-|x-y|$ around $|x-y|.$ To that end, we need boundedness of $v_j,$ which we do not have apriorily, but we can achieve it through truncation. Denoting
\begin{equation}
\label{3.21} 
f_j(x)=\frac{1}{a_j}\left(\int_{Q} \left|\Go(u_j)(x,y)\right|^pdy\right)^{1/p},
\end{equation}
we have that 
\begin{equation}
\label{3.22} 
\|f_j\|_{L^p(Q)}\longrightarrow 0
\end{equation}
by the initial assumption. For almost every pair $x,y\in Q$ we have by the triangle inequality
\begin{equation}
\label{3.23} 
\bigl||u_j(x)-y|-|x-y| \bigr| \le\Go(u_j)(x,y)+|u_j(y)-y|.
\end{equation}
Next we apply the following averaged reverse triangle inequality, the proof of which is routine and left to the reader. 
\textit{For every $x\in Q$ and 
$z\in\mathbb R^n$ one has}
\begin{equation}
\label{3.24}
 |z-x|\le C_0(n)\left(\int_{Q} \bigl||z-y|-|x-y|\bigr|^p dy\right)^{1/p}.
\end{equation}
Taking now $z=\Bu_j(x)$ in (\ref{3.21})-(\ref{3.22}) and using Minkowski's inequality in $y,$ we get
\begin{equation}
\label{3.25}
 |v_j(x)|\le C_0(n)\bigl(f_j(x)+1\bigr)\qquad\text{for a.e.} \qquad x\in Q.
\end{equation}
In the next step we replace $v_j$ and $u_j$ by equivalent $L^\infty$ sequences by truncation as follows. Let $w_j$ 
be the radial truncation of $v_j$ at the cutoff $2C_0(n)$ and set $\widetilde u_j=x+a_jw_j.$ By (\ref{3.23}), 
we have the obvious bounds 
\begin{equation}
\label{3.26}
 \|w_j\|_{L^\infty(Q)}\le2 C_0(n) \qquad\text{and}\qquad |v_j-w_j|=(|v_j|-2C_0(n))_{+}\le C_0(n)f_j,
\end{equation}
thus from (\ref{3.22}) we get
\begin{equation}
\label{3.27}
\|v_j-w_j\|_{L^p(Q)}\longrightarrow  0.
\end{equation} 
Next, we have by the Minkowski inequality
\begin{equation}
\label{3.28}
\|E w_j-E v_j\|_{L^p(Q\times Q)}\le C(n,p)\|v_j-w_j\|_{L^p(Q)}\longrightarrow  0.
\end{equation}
The reverse triangle inequality implies
\begin{equation}
\label{3.29}
\frac{1}{a_j} \|\Go(u_j)-\Go(\widetilde u_j)\|_{L^p(Q\times Q)}\le2\|v_j-w_j\|_{L^p(Q)}\longrightarrow  0. 
\end{equation}
We can write a Taylor expansion inequality for $|x|$ as follows:
\begin{equation}
\label{3.30}
0\leq |z+h|-|z|-\frac{z\cdot h}{|z|} \leq \min\left(2|h|,\frac{|h|^2}{2|z|}\right).
\end{equation}
Thus, thanks to $|w_j(x)-w_j(y)|\le 4C_0(n),$ implementing (\ref{3.30}) on $z=x-y,$ $h=a_j(w_j(x)-w_j(y))$ we obtain by the 
dominated convergence theorem and the Minkowski inequality, that
\begin{equation}
\label{3.31}
 \frac{1}{a_j}\|\Go(\widetilde u_j)-a_jE w_j\|_{L^p(Q\times Q)} \longrightarrow 0.
\end{equation}
A combination of (\ref{3.28}), (\ref{3.29}), and (\ref{3.31}) gives 
\begin{equation}
\label{3.32}
 \frac{1}{a_j}\|\Go(u_j)-a_jE v_j\|_{L^p(Q\times Q)} \longrightarrow 0,
 \end{equation}
which, taking into account (\ref{3.14}), yields the strong convergence 
\begin{equation}
\label{3.33}
 \|E v_j\|_{L^p(Q\times Q)} \longrightarrow 0.
 \end{equation}
 Now, owing back to the unweighted Korn's first inequality on cubes in (\ref{3.12}), we get a sequence of skew symmetric 
 matrices $\BA_j$ and translations $c_j\in\mathbb R^n$ with
\begin{equation}
\label{3.34}
 \|v_j-\BA_jx-c_j\|_{L^p(Q)} \longrightarrow 0.
\end{equation}
Their coefficients are bounded, i.e., $\|\BA_j\|_{L^\infty}\leq M,$ because $\|v_j\|_{L^p}=1$ and the triangle inequality and
 the fact, that norms on a finite-dimensional space of affine maps are equivalent. Now we can undo the linearization in (\ref{3.14})
 by combining $x$ and $a_j\BA_jx$ to get a better approximation of $u_j.$ The mapping $\BI+a_j\BA_j$ has to be the first order approximation of a proper rotation, and proper rotations are exponentials of skew-symmetric matrixes, thus due to Taylor expansion
 of exponentials, we have to define the better approximation as 
\begin{equation}
\label{3.35}
 \BB_j=\exp(a_j\BA_j)\in\mathrm{SO}(n).
\end{equation}
Thanks to the uniform bound $\|\BA_j\|_{L^\infty}\leq M,$ the exponential satisfies 
$$
\|\BB_j-\BI-a_j\BA_j\|\le C(M,n)a_j^2
$$ 
by Taylor expansion and the fact that $a_j\to 0.$ Hence we can estimate
\begin{align*}
 \|u_j-\BB_jx-a_jc_j\|_{L^p(Q)} &\le a_j\|v_j-\BA_jx-c_j\|_{L^p(Q)} +\|(\BI+a_j\BA_j-\BB_j)x\|_{L^p(Q)}\\
 &=o(a_j).
\end{align*}
This contradicts the fact that the distance of $u_j$ from all affine isometries is $a_j$. This concludes the proof of Theorem~2.5.
\end{proof}

\vspace{0.2cm}

\begin{remark}
\label{Rem:3.2}
To achieve full analogy with the non-fractional case one would wish to replace the left hand side of (\ref{2.11}) 
with $\inf_{\BR\in SO(n)}[u-\BR x]_{W^{s,p}(\GO)}$. To achieve this one should restrict 
the space $W^{s,p}$ and a natural candidate would be the subset of the $W^{s,p}$-closure of positively oriented deformations: \begin{equation*} 
 W^{s,p}_+(\Omega;\mathbb R^n)
 :=\overline{\left\{v\in C^\infty(\Omega;\mathbb R^n)
 \cap W^{s,p}(\Omega;\mathbb R^n):
 \det Dv>0\ \hbox{in }\Omega\right\}}^{\,W^{s,p}}
\end{equation*} 
Is the inequality 
\begin{equation}
\label{3.36}
 \inf_{\BR \in \mathrm{SO}(n)}[u-\BR x]_{W^{s,p}(\GO)}^p \leq C \int_{\GO}\int_{\GO} \frac{\left| \Go(u)\right|^p}{|x-y|^{n+ps}}dxdy. 
\end{equation}
true for $u\in W^{s,p}_+(\GO; \mathbb R^n)$? 

It is not hard to see that (\ref{3.36}) does not hold for $sp<1$ (because, by the weak topology, any reflection belongs to $W^{s,p}_+(\GO ; \mathbb R^n)$), and that it does hold for $u\in W^{\frac{n-1}{n},n}$ (due to the continuous extension of the Jacobian; see [\ref{bib:Bre.Ngu.}]). The authors will address this question in more detail in a 
forthcoming work.
\end{remark}


\section*{Acknowledgments}
D. H.'s work is partially supported by the National Science Foundation under Grants No. DMS-2206239. 
We would like to thank Tadele Mengesha and James Scott for many interesting discussions on the subject. 
We used ChatGpt Astra-6 for reference search and for some help with the strong convergence of $Ev_j\to 0$ in 
Theorem~2.3. The approach in (\ref{3.21}) and (\ref{3.24}) was suggested by AI with full human assistance on 
which direction to go.



\begin{thebibliography}{999}





\bibitem{} G. Acosta, I. Drelichman, R. Dur\'an, F. L\'opez-Garcia, and I. Ojea. The fractional Korn inequality on uniform domains and new Korn inequalities for truncated seminorms, \textit{preprint}, https://arxiv.org/pdf/2601.08096
\label{bib:Aco.Dre.Dur.Lop.Oje.}


\bibitem{} Bourgain, J., Brezis, H., Mironescu, P., Another look at Sobolev spaces, \textit{Optimal Control and Partial Differential Equations. } A Volume in Honor of A. Bensoussan's 60th Birthday, IOS Press, Amsterdam, 2001, pp. 439--455.
\label{bib:Bou.Bre.Mir.1}

\bibitem{}  H. Brezis and H.-M. Nguyen. The Jacobian determinant revisited, \textit{Inventiones Mathematicae,} 
185 (2011), 17--54. 
\label{bib:Bre.Ngu.}

\bibitem{} M. Chermisi and S. Conti. Multiwell Rigidity in Nonlinear Elasticity, \textit{SIAM J. Math. Anal.,} Vol. 42, Iss. 5, (2010).
\label{bib:Che.Con.}


\bibitem{}  S. Conti, G. Dolzmann and S. M\"uller. Korn's second inequality and geometric rigidity with mixed growth conditions. \textit{Calc. Var. PDE,} (2014) 50:437--454.
\label{bib:Con.Dol.Mue.}


\bibitem{} S. Conti and B. Zwicknagl.  The tapering length of needles in martensite/martensite macrotwins, \textit{Archive for Rational Mechanics and Analysis,} 247:63 (2023).
\label{bib:Con.Zwi.}


\bibitem{} C. De Lellis a, L\'aszl\'o Szekelyhidi Jr. Simple proof of two-well rigidity, \textit{C. R. Acad. Sci. Paris, Ser. I} 
343 (2006) 367--370.
\label{bib:DeL.Sze.}


\bibitem{}  I. Drelichman and R. G. Dur\'an. Improved Poincaré inequalities with weights,
\textit{Journal of Mathematical Analysis and Applications,} 347 (2008), no. 1, 286--293.
\label{bib:Dre.Dur.}



\bibitem{} E. Di Nezza, G. Palatucci, and E. Valdinoci. Hitchhiker's guide to the fractional Sobolev spaces, \textit{Bulletin Des Sciences Math\'ematiques,} 
Vol. 136, Issue 5, July--August 2012, pp. 521--573. 
\label{bib:DiN.Pat.Val.}



\bibitem{} Q. Du and K. Zhou, Mathematical analysis for the peridynamic non-local continuum theory. \textit{ESAIM: Mathematical Modeling and Numerical Analysis,} Vol. 45 (2011) no. 2, pp. 217--234.
\label{bib:Du.Zho.}


\bibitem{} L.C. Evans and R. F. Gariepy. \textit{Measure theory and fine properties of functions,} Revised Edition, CRC Press
Taylor and Francis Group, Boka Raton \  London \ New York. 
\label{bib:Eva.Gar.}


\bibitem{}  M. Friedrich, L. Kreutz, and K. Zemas. Geometric rigidity in variable domains and derivation of linearized models for elastic materials with free surfaces, \textit{Annales de l'Institut Henri Poincar\'e C,} 42 : 5, pp: 1093--1163
\label{bib:Fri.Kre.Zem.} 


\bibitem{} G. Friesecke, R. D. James, and S. M\"uller. A theorem on geometric rigidity and the derivation of nonlinear plate theory
 from three-dimensional elasticity. \textit{Comm. Pure Appl. Math.}, 55(11):1461-1506, 2002.
\label{bib:Fri.Jam.Mue.}


\bibitem{} F. W. Gehring and B. G. Osgood. Uniform domains and the quasi-hyperbolic metric, \textit{J. Analyse Math.}
 36 (1979), 50--74.
\label{bib:Geh.Osg.}


\bibitem{} H. Hanche-Olsen, H. Holden. The Kolmogorov--Riesz compactness theorem,\textit{Expositiones Mathematicae,}
Vol. 28, Iss. 4, pp: 385--394, 2010.
\label{bib:Han.Hol.} 



\bibitem{} D. Harutyunyan, T. Mengesha, H. Mikayelyan, and J. M. Scott. Fractional Korn's inequalities without boundary conditions. \textit{Mathematics and Mechanics of Complex Systems,} 11(4):497--524, 2023.
\label{bib:Har.Men.Mik.Sco.}



\bibitem{} D. Harutyunyan and H. Mikayelyan. On the fractional Korn inequality in bounded domains: Counterexamples to the case 
$ps<1$. \textit{Advances in Nonlinear Analysis,} 12(1):20220283, 2023.
\label{bib:Har.Mik.}


\bibitem{} D. Harutyunyan and A. M. Rodrigues. On weighted rigidity estimates for bulk and thin domains, \textit{preprint}, 
https://arxiv.org/abs/2608.00863
\label{bib:Har.Rod.}


\bibitem{} D. Harutyunyan. Singularity removal in nonlocal operator inequalities on domains and applications, \textit{preprint}, 
https://arxiv.org/abs/2610.39437
\label{bib:Harutyunyan}

\bibitem{}  R. Hurri--Syrj\"anen. An Improved Poincar\'e Inequality, \textit{ Proc. Amer. Math. Soc.,} 120 (1), pp. 213--222, 1994.
\label{bib:Hurri}


\bibitem{}  F. John. Rotation and strain, \textit{Comm. Pure Appl. Math.} 14 (1961), 391--413.
\label{bib:John}

\bibitem{}  P. W. Jones. Quasiconformal mappings and extendability of functions in Sobolev spaces. \textit{Acta Math.,} 1
47:71--88, 1981.
\label{bib:Jones}


\bibitem{} O. Marito and J. Sarvas. Injectivity theorems in plane and space, \textit{Annales Academire Scientiarum Fennicre Series A. I. Mathematica} Volumen 4, 1978/1979, 383--401.
\label{bib:Mar.Sar.}

\bibitem{}  V. Maz'ya and T. Shaposhnikova. On the Bourgain, Brezis, and Mironescu theorem concerning limiting embeddings of fractional Sobolev spaces, \textit{J. Funct. Anal.} 195 (2) (2002) 230--238; Erratum, \textit{J. Funct. Anal. }201 (1) (2003) 298--300.
\label{bib:Maz.Sha.}

\bibitem{}  T. Mengesha. Fractional Korn and Hardy-type inequalities for vector fields in half space, \textit{Communications in Contemporary Mathematics,} 
Vol. 21, No. 7 (2019)
\label{bib:Mengesha.1}

\bibitem{}  T. Mengesha and J. M. Scott. A Fractional Korn-type inequality for smooth domains and a regularity estimate for nonlinear nonlocal systems of equations,  \textit{Commun. Math. Sci.} Vol. 20, N0. 2, 405--423, 2022.
\label{bib:Men.Sco.1}


\bibitem{}  T. Mengesha and J. M. Scott. Consistency in fractional Korn's inequalities on Lipschitz domains, \textit{preprint.}
https://arxiv.org/abs/2610.33325
\label{bib:Men.Sco.2}


\bibitem{} M\"uller, S, Scardia, L, and  Zeppieri, C. I. Geometric rigidity for incompatible fields, and an application to strain-gradient plasticity, \textit{Indiana University Mathematics Journal,} vol. 63, no. 5, pp. 1365--1396.
\label{bib:Mue.Sca.Zep.}


\bibitem{} A. Rutkowski. Fractional Korn's inequality on subsets of the Euclidean space,
\textit{Mathematical Inequalities And Applications,} Vol. 25, Nr. 2 (2022), 359--367.
\label{bib:Rutkowski.1}


\bibitem{}  J. M. Scott. Nonlocal geometric rigidity inequalities and the justification of linearized peridynamics, \textit{preprint}
https://arxiv.org/abs/2610.28268
\label{bib:Scott}

\bibitem{} A. Seeger. A note on Triebel-Lizorkin spaces. \textit{Banach Center Publications,} 22(1):391--400, 1989.
\label{bib:Seeger}

 
\bibitem{}  S. A. Silling. Reformulation of elasticity theory for discontinuities and long-range forces, \textit{J. Mech. Phys. Solids,} 48, 175--209.
\label{bib:Silling.1}

\bibitem{}  S. A. Silling. Linearized theory of peridynamic states, \textit{J. Elastisity,} 99, 8--111.
\label{bib:Silling.2}

\bibitem{}  S. A. Silling, M. Epton, O. Weckner, J. Xu, and E. Askari. Peridynamic states and constitutive modeling, \textit{J. Elastisity,} 88, 151--184.
\label{bib:Sil.Ept.Wec.Xu.Ask.}


\bibitem{} R. Temam and A. Miranville, \textit{Mathematical Modeling in Continuum Mechanics,} Cambridge University Press, 2001. 
\label{bib:Tem.Mir.}


\bibitem{} J. Zaliaduonis and S. Gatidis. A Probabilistic Generalization of the Mazur-Ulam Theorem, \textit{Preprint,} 
https://arxiv.org/pdf/2601.03900
\label{bib:Zal.Gat.}

\end{thebibliography}
\end{document}